\RequirePackage{ifpdf}
\ifpdf 
\documentclass[pdftex]{sigma}
\else
\documentclass{sigma}
\fi

\begin{document}

\allowdisplaybreaks

\renewcommand{\PaperNumber}{095}

\FirstPageHeading

\ShortArticleName{Roots of the Yablonskii--Vorob'ev Polynomials}

\ArticleName{Irrationality of the Roots of the Yablonskii--Vorob'ev Polynomials and Relations between Them}

\Author{Pieter ROFFELSEN}

\AuthorNameForHeading{P. Rof\/felsen}

\Address{Radboud Universiteit Nijmegen, IMAPP, FNWI,\\ Heyendaalseweg 135, 6525 AJ Nijmegen, the Netherlands}
\Email{\href{mailto:roffelse@science.ru.nl}{roffelse@science.ru.nl}}

\ArticleDates{Received November 13, 2010, in f\/inal form December 08, 2010;  Published online December 14, 2010}

\Abstract{We study the Yablonskii--Vorob'ev polynomials, which are special polynomials used to represent rational solutions of the second Painlev\'e equation. Divisibility properties of the coef\/f\/icients of these polynomials, concerning powers of 4, are obtained and we prove that the nonzero roots of the Yablonskii--Vorob'ev polynomials are irrational. Furthermore, relations between the roots of these polynomials for consecutive degree are found by considering power series expansions of rational solutions of the second Painlev\'e equation.}

\Keywords{second Painlev\'e equation; rational solutions; power series expansion; irrational roots; Yablonskii--Vorob'ev polynomials}

\Classification{34M55}

\section{Introduction}

In this paper we study the Yablonskii--Vorob'ev polynomials $Q_n$, with a special interest in their roots. These polynomials were derived by Yablonskii and Vorob'ev, while examining the hierarchy of rational solutions of the second Painlev\'e equation. The Yablonskii--Vorob'ev polynomials are def\/ined by the dif\/ferential-dif\/ference equation
\begin{gather} \label{eqqn}
Q_{n+1} Q_{n-1}=zQ_n^2-4(Q_n Q_n''-(Q_n')^2),\end{gather}
with $Q_0=1$ and $Q_1=z$. From the recurrence relation, it is clear that the functions $Q_n$ are rational, though it is far  from obvious that they are polynomials, since in every iteration one divides by $Q_{n-1}$. The Yablonskii--Vorob'ev polynomials $Q_n$ are monic polynomials of degree $\frac{1}{2}n(n+1)$, with integer coef\/f\/icients. The f\/irst few are given in Table~\ref{tableQn}.
\begin{table}[t]
\centering \caption{}\label{tableQn} \vspace{1mm}

\begin{tabular}{r@{\,}c@{\,}l }
\hline
\multicolumn{3}{c}{Yablonskii--Vorob'ev polynomials\tsep{1pt}\bsep{1pt}}\\
\hline
$Q_2$ & $=$ & $4+z^3$\tsep{1pt}\\
$Q_3$ & $=$ & $-80+20 z^3+z^6$\\
$Q_4$ & $=$ & $z(11200+60 z^6+z^9)$\\
$Q_5$ & $=$ & $-6272000-3136000 z^3+78400 z^6+2800 z^9+140 z^{12}+z^{15}$\\
$Q_6$ & $=$ & $-38635520000+19317760000 z^3+1448832000 z^6-17248000 z^9+
627200z^{12}$\\
 & & $+18480 z^{15}+280 z^{18}+z^{21}$\\
$Q_7$ & $=$ & $z (-3093932441600000-49723914240000 z^6-828731904000 z^9+13039488000
z^{12}$\\
 & & $+62092800 z^{15}+5174400 z^{18}+75600 z^{21}+504 z^{24}+z^{27})$\\
$Q_8$ & $=$ & $-991048439693312000000-743286329769984000000 z^3$\\
 & & $+37164316488499200000
z^6+1769729356595200000 z^9+126696533483520000 z^{12}$\\
 & & $+407736096768000
z^{15}-6629855232000 z^{18}+124309785600 z^{21}+2018016000 z^{24}$\\
 & & $+32771200z^{27}+240240 z^{30}+840 z^{33}+z^{36}$\bsep{1pt}\\ \hline
\end{tabular}

\end{table}

Yablonskii \cite{yablonskii} and Vorob'ev \cite{vorob} expressed the rational solutions of the second Painlev\'e equation,
\[P_{\rm II}(\alpha):\ \ w''(z)=2w(z)^3+zw(z)+\alpha,\]
with complex parameter $\alpha$, in terms of the Yablonskii--Vorob'ev polynomials, as summerized in the following theorem:

\begin{theorem}
\label{thmYV}
$P_{\rm II}(\alpha)$ has a rational solution iff $\alpha=n\in\mathbb{Z}$. For $n\in\mathbb{Z}$ the rational solution is unique and if $n\geq 1$, then it is equal to
\[w_n=\frac{Q_{n-1}'}{Q_{n-1}}-\frac{Q_n'}{Q_n}.\]
The other rational solutions are given by $w_0=0$ and for $n\geq 1$, $w_{-n}=-w_n$.
\end{theorem}

The rational solutions of $P_{\rm II}$ can also be determined, using the B\"acklund transformations, f\/irst given by Gambier~\cite{gambier}, of the second Painlev\'e equation, by
\begin{gather*}
w_{n+1}= -w_n-\frac{2n+1}{2w_n^2+2w_n'+z},\qquad
w_{-n}= -w_n,
\end{gather*}
with ``seed solution'' $w_0=0$; see also Lukashevich~\cite{lukashevich} and Noumi~\cite{noumi}.

We note that the Yablonskii--Vorob'ev polynomials f\/ind many applications in physics. For instance, solutions of the Korteweg--de Vries equation (Airault, McKean and Moser~\cite{airault}) and the Boussinesq equation (Clarkson~\cite{clarksonboussinesq}) can be expressed in terms of these polynomials. Clarkson and Mansf\/ield~\cite{clarksonmansfield} studied the structure of the roots of the Yablonskii--Vorob'ev polynomials~$Q_n$ and observed that the roots, of each of these polynomials, form a highly regular triangular-like pattern, for~$n\leq 7$, suggesting that they have interesting properties. This further motivates studying the zeros of the Yablonskii--Vorob'ev polynomials.

In Section~\ref{section2} the divisibility of the coef\/f\/icients of the Yablonskii--Vorob'ev polynomials by powers of~$4$ is examined. From the divisibility properties found, we conclude that nonzero roots of the Yablonskii--Vorob'ev polynomials are irrational. In Section~\ref{section3} we study power series expansions of (functions related to) the rational solution $w_n$ of $P_{\rm II}(n)$, around poles of $w_n$. This leads to relations between the roots of $Q_{n-1}$ and $Q_n$. These relations suggest deeper connections between the zeros of $Q_{n-1}$ and $Q_n$. Similarly, we look at power series expansions of (functions related to) the rational solution $w_n$ of $P_{\rm II}(n)$ around $0$, in Section~\ref{section4}. We obtain polynomial expressions in~$n$, with rational coef\/f\/icients, for sums of f\/ixed negative powers of the nonzero roots of~$Q_n$.

\section{Nonzero roots are irrational}\label{section2}

The Yablonskii--Vorob'ev polynomials $Q_n$ are monic polynomials of degree $\frac{1}{2}n(n+1)$, and Taneda~\cite{taneda} proved:
\begin{itemize}\itemsep=0pt
\item if $n\equiv 1 \pmod{3}$, then $\frac{Q_n}{z}\in\mathbb{Z}[z^3]$;
\item if $n\not\equiv 1 \pmod{3}$, then $Q_n\in\mathbb{Z}[z^3]$.
\end{itemize}
Therefore, we have{\samepage
\begin{gather}
\label{eqcoef}
Q_n=z^{\frac{1}{2}n(n+1)}+a_1^nz^{\frac{1}{2}n(n+1)-3}+a_2^nz^{\frac{1}{2}n(n+1)-6}+\cdots+a_{\left[\frac{1}{6}n(n+1)\right]}^nz^{\frac{1}{2}n(n+1)-3\left[\frac{1}{6}n(n+1)\right]},
\end{gather}
for certain $a_s^n\in\mathbb{Z}$, with convention $a_0^n=1$, where $\left[ \cdot \right]$ denotes the f\/loor function.}

\begin{lemma}
\label{coeflemma}
For every $0\leq m\leq \left[\frac{1}{6}n(n+1)\right]$, we have $4^m\mid a_m^n$.
\end{lemma}
\begin{proof}
We proceed by proving the following statement, by induction, for all $M\in\mathbb{N}$:

For every $1\leq m\leq M$, for all $n\in\mathbb{N}$, whenever $m\leq \left[\frac{1}{6}n(n+1)\right]$, we have $4^m\mid a_m^n$, and
\begin{gather*}
4^M\mid a_{M+1}^n ,\qquad  4^M\mid a_{M+2}^n ,\qquad   \ldots, \qquad
 4^M\mid a_{\left[\frac{1}{6}n(n+1)\right]}^n.
\end{gather*}
Observe that the case $M=0$ is trivial. Now suppose the statement is true for $M\in\mathbb{N}$.
Then there are $b_s^n\in\mathbb{Z}$, such that for every $n\in\mathbb{N}$,
\[Q_n=z^{\frac{1}{2}n(n+1)}+4b_1^nz^{\frac{1}{2}n(n+1)-3}+4^2b_2^nz^{\frac{1}{2}n(n+1)-6}+\cdots+4^Mb_M^n z^{\frac{1}{2}n(n+1)-3M}+4^M P_n,\]
where $P_n\in\mathbb{Z}[z]$ is zero or has degree less or equal to $\frac{1}{2}n(n+1)-3(M+1)$, and if $m>\left[\frac{1}{6}n(n+1)\right]$, then $b_m^n=0$.

To complete the induction, we need to show that for every $n\in\mathbb{N}$, $4\mid P_n$. We prove this by induction with respect to $n$. Observe that $P_0=0$ and $P_1=0$, therefore, indeed $4\mid P_0$ and $4\mid P_1$.
Assume $4\mid P_{n-1}$ and $4\mid P_n$. Then $4^MP_n\equiv 0 \pmod{4^{M+1}}$, therefore, modulo $4^{M+1}$, we have:
\begin{gather*}
z^{\max(0,n(n+1)-3M+1)}\mid zQ_n^2,\qquad
z^{\max(0,n(n+1)-3M+1)}\mid 4Q_nQ_n'',\\
z^{\max(0,n(n+1)-3M+1)}\mid 4(Q_n')^2.
\end{gather*}
By the def\/inition of $Q_{n+1}$ \eqref{eqqn},
\[
Q_{n+1} Q_{n-1}=zQ_n^2-4\big(Q_n Q_n''-(Q_n')^2\big),
\]
so \begin{gather}
\label{eqdiv}
z^{\max(0,n(n+1)-3M+1)}\mid Q_{n+1}Q_{n-1} \pmod{4^{M+1}}.
\end{gather}

Let us consider $Q_{n+1}Q_{n-1}$. Since $4\mid P_{n-1}$, we have
\[
4^MP_{n-1}\equiv 0 \pmod{4^{M+1}},
\] therefore, modulo $4^{M+1}$,
\begin{gather}
Q_{n+1}Q_{n-1}\equiv Q_{n+1}z^{\frac{1}{2}n(n-1)}+Q_{n+1}\big(4b_1^{n-1}z^{\frac{1}{2}n(n-1)-3}\nonumber\\
\phantom{Q_{n+1}Q_{n-1}\equiv}{} +4^2b_2^{n-1}z^{\frac{1}{2}n(n-1)-6}+\cdots+4^Mb_M^{n-1} z^{\frac{1}{2}n(n-1)-3M}\big).\label{eqdiv2}
\end{gather}
Since
\begin{gather*}
Q_{n+1}= z^{\frac{1}{2}(n+1)(n+2)}+4b_1^{n+1}z^{\frac{1}{2}(n+1)(n+2)-3}\\
 \phantom{Q_{n+1}=}{} +4^2b_2^{n+1}z^{\frac{1}{2}(n+1)(n+2)-6}+\cdots+4^Mb_M^{n+1} z^{\frac{1}{2}(n+1)(n+2)-3M}+4^M P_{n+1},
\end{gather*}
we have, modulo $4^{M+1}$, \begin{gather*}
z^{\max(0,n(n+1)-3M+1)}\mid Q_{n+1}\big(4b_1^{n-1}z^{\frac{1}{2}n(n-1)-3}+4^2b_2^{n-1}z^{\frac{1}{2}n(n-1)-6}\\
\qquad{} +\cdots+4^Mb_M^{n-1} z^{\frac{1}{2}n(n-1)-3M}\big).
\end{gather*}
Hence, by \eqref{eqdiv} and \eqref{eqdiv2}, \[
z^{\max(0,n(n+1)-3M+1)}\mid Q_{n+1}z^{\frac{1}{2}n(n-1)} \pmod{4^{M+1}},
\]
which implies
\[
z^{\max(0,\frac{1}{2}(n+1)(n+2)-3M)} \mid Q_{n+1} \pmod{4^{M+1}}.\]
Since
\begin{gather*}
Q_{n+1}= z^{\frac{1}{2}(n+1)(n+2)}+4b_1^{n+1}z^{\frac{1}{2}(n+1)(n+2)-3}\\
\phantom{Q_{n+1}=}{} +4^2b_2^{n+1}z^{\frac{1}{2}(n+1)(n+2)-6}+\cdots+4^Mb_M^{n+1} z^{\frac{1}{2}(n+1)(n+2)-3M}+4^M P_{n+1},
\end{gather*}
we have, therefore, $4\mid P_{n+1}$. Hence, by induction, for all $n\in\mathbb{N}$, $4\mid P_n$.

The lemma follows by induction on $M$.
\end{proof}

Let us denote the coef\/f\/icient of the lowest degree term in $Q_n$ by \[
x_n:=a_{\left[\frac{1}{6}n(n+1)\right]}^n,
\] i.e.~$x_n$ is the constant coef\/f\/icient in $Q_n$ if $n\not\equiv 1 \pmod{3}$, and $x_n$ is the coef\/f\/icient of $z$ in $Q_n$ if $n\equiv 1 \pmod{3}$.
Fukutani, Okamoto, and Umemura~\cite{fukutani} proved that the roots of the Yablonskii--Vorob'ev polynomials are simple, hence~$x_n$ is nonzero.
Let $p_n$ be the multiplicity of~2 in the prime factorization of $x_n$.
As a consequence of Lemma~\ref{coeflemma}, we obtain that $p_n\geq 2\left[\frac{1}{6}n(n+1)\right]$. We prove \[p_n=\left[\frac{1}{3}n(n+1)\right].\]

Observe that $x_n=Q_n(0)$ if $n\not\equiv 1\pmod{3}$, and $x_n=Q_n'(0)$ if $n\equiv 1\pmod{3}$.
Fukutani, Okamoto, and Umemura \cite{fukutani} derived the following identity for the Yablonskii--Vorob'ev polynomials:
\[Q_{n+1}' Q_{n-1}-Q_{n+1}Q_{n-1}' = (2n+1)Q_n^2.\] Using this identity at $0$, we obtain
\[
x_{n+1}x_{n-1}=\begin{cases}(2n+1)x_n^2 &\text{ if $n\equiv 0 \pmod{3}$,}\\
-(2n+1)x_n^2 &\text{ if $n\equiv 2 \pmod{3}$.} \end{cases}
\]
By evaluating equation~\eqref{eqqn} at $0$,
\[x_{n+1}x_{n-1}=4x_n^2,\qquad  \text{if $n\equiv 1 \pmod{3}$.}\]
Therefore, we have the following recursion for $(x_n)_n$:
\begin{gather*}
 x_0=1 ,
 \qquad x_1=1 \qquad{\rm and}\\
 x_{n+1}x_{n-1}=
\begin{cases}(2n+1)x_n^2 &\text{ if $n\equiv 0 \pmod{3}$,}\\
4x_n^2 &\text{ if $n\equiv 1 \pmod{3}$,}\\
-(2n+1)x_n^2 &\text{ if $n\equiv 2 \pmod{3}$.} \end{cases}
\end{gather*}
So, we obtain the following recursion for $(p_n)_n$:
\begin{gather*}
 p_0=0 ,  \qquad p_1=0 \qquad {\rm  and}\\
 p_{n+1}=\begin{cases} 2p_n-p_{n-1} & \text{if $n\not\equiv 1 \pmod{3}$,}\\
2+2p_n-p_{n-1} & \text{if $ n\equiv 1 \pmod{3}$.}
\end{cases}
\end{gather*}
Using this recursion, the formula $p_n=\left[\frac{1}{3}n(n+1)\right]$, can be proven directly, by induction.

\begin{remark}
Kaneko and Ochiai \cite{kaneko} found an explicit expression for the coef\/f\/icients $x_n$. But deriving the formula $p_n=\left[\frac{1}{3}n(n+1)\right]$ directly from this expression seems to be a dif\/f\/icult task.
\end{remark}

\begin{theorem} The nonzero roots of the Yablonskii--Vorob'ev polynomials are irrational.
\end{theorem}

\begin{proof}
Let $n\not\equiv 1 \pmod{3}$.
Suppose $x$ is a rational root of $Q_n$. Since $Q_n\in\mathbb{Z}[z]$ is monic, by Gauss's lemma, $x\in\mathbb{Z}$.
By Lemma~\ref{coeflemma}, \[Q_n\equiv z^{\frac{1}{2}n(n+1)} \pmod{4},\] so $x$ is even. Let $y:=\frac{x}{2}$, then, by equation~\eqref{eqcoef},
\[0=(2y)^{\frac{1}{2}n(n+1)}+a_1^n(2y)^{\frac{1}{2}n(n+1)-3}+a_2^n(2y)^{\frac{1}{2}n(n+1)-6}+\cdots+a_{\frac{1}{6}n(n+1)-1}^n(2y)^3+a_{\frac{1}{6}n(n+1)}^n.\]
By Lemma \ref{coeflemma}, for every $m\leq \frac{1}{6}n(n+1)$, we have $4^m\mid a_m^n$.
Hence
\begin{gather*}
2^{\frac{1}{2}n(n+1)}\mid(2y)^{\frac{1}{2}n(n+1)} ,\qquad
 2^{\frac{1}{2}n(n+1)-1}\mid a_1^n(2y)^{\frac{1}{2}n(n+1)-3} ,\\
 2^{\frac{1}{2}n(n+1)-2}\mid a_2^n(2y)^{\frac{1}{2}n(n+1)-6} ,\qquad
\ldots, \qquad 2^{\frac{1}{2}n(n+1)-\frac{1}{6}n(n+1)+1}\mid a_{\frac{1}{6}n(n+1)-1}(2y)^3.
\end{gather*}
 So \[2^{\frac{1}{3}n(n+1)+1}\mid a_{\frac{1}{6}n(n+1)}^n=x_n,\]
which implies \[p_n\geq \frac{1}{3}n(n+1)+1.\]
But $p_n=\frac{1}{3}n(n+1)$, a contradiction, hence roots of $Q_n$ are irrational.

If $n\equiv 1 \pmod{3}$, we can apply the same reasoning to $\frac{Q_n}{z}$, and show that roots of $\frac{Q_n}{z}$ are irrational. Therefore, nonzero roots of $Q_n$ are irrational.
\end{proof}

This result raises the question whether the Yablonskii--Vorob'ev polynomials, excluding the trivial factor $z$ in case $n\equiv 1 \pmod{3}$, are irreducible in $\mathbb{Q}[z]$. Kametaka \cite{kametaka} showed that for $n\leq 23$, the Yablonskii--Vorob'ev polynomials $Q_n$ are indeed irreducible.

\section[Relations between roots of the Yablonskii-Vorob'ev polynomials]{Relations between roots of the Yablonskii--Vorob'ev\\ polynomials}\label{section3}

By Theorem \ref{thmYV}, for $n\geq 1$, the unique rational solution of $P_{\rm II}(n)$ is given by
\[w_n=\frac{Q_{n-1}'}{Q_{n-1}}-\frac{Q_n'}{Q_n}.\]
Fukutani, Okamoto, and Umemura \cite{fukutani} proved that the roots of the Yablonskii--Vorob'ev polynomials are simple, hence
\begin{gather}
\label{eqyn}
w_n= \sum_{k=1}^{\frac{1}{2}n(n-1)}{\frac{1}{z-z_{n-1,k}}}- \sum_{k=1}^{\frac{1}{2}n(n+1)}{\frac{1}{z-z_{n,k}}},
\end{gather}
where the $z_{m,k}$ are the roots of $Q_m$. From equation~\eqref{eqyn} and the fact that $w_n$ is the rational solution of $P_{\rm II}(n)$, we obtain relations between the zeros of $Q_{n-1}$ and $Q_n$.

\begin{theorem}
\label{TREL}
For $1\leq j\leq \frac{1}{2}n(n-1)$:
\begin{gather*}
\sum_{k=1,\; k\neq j}^{\frac{1}{2}n(n-1)}{\frac{1}{z_{n-1,j}-z_{n-1,k}}}- \sum_{k=1}^{\frac{1}{2}n(n+1)}{\frac{1}{z_{n-1,j}-z_{n,k}}}=0,\\
 \sum_{k=1,\;  k\neq j}^{\frac{1}{2}n(n-1)}{\frac{1}{(z_{n-1,j}-z_{n-1,k})^2}}- \sum_{k=1}^{\frac{1}{2}n(n+1)}{\frac{1}{(z_{n-1,j}-z_{n,k})^2}}=\frac{z_{n-1,j}}{6},\\
 \sum_{k=1,\; k\neq j}^{\frac{1}{2}n(n-1)}{\frac{1}{(z_{n-1,j}-z_{n-1,k})^3}}-
 \sum_{k=1}^{\frac{1}{2}n(n+1)}{\frac{1}{(z_{n-1,j}-z_{n,k})^3}}=-\frac{n+1}{4},\\
 \sum_{k=1,\; k\neq j}^{\frac{1}{2}n(n-1)}{\frac{1}{(z_{n-1,j}-z_{n-1,k})^5}}- \sum_{k=1}^{\frac{1}{2}n(n+1)}{\frac{1}{(z_{n-1,j}-z_{n,k})^5}}=z_{n-1,j}\left(\frac{n+1}{24}-\frac{1}{36}\right).
\end{gather*}
For $1\leq j\leq \frac{1}{2}n(n+1)$:
\begin{gather*}
\sum_{k=1}^{\frac{1}{2}n(n-1)}{\frac{1}{z_{n,j}-z_{n-1,k}}}- \sum_{k=1,\; k\neq j}^{\frac{1}{2}n(n+1)}{\frac{1}{z_{n,j}-z_{n,k}}}=0,\\
 \sum_{k=1}^{\frac{1}{2}n(n-1)}{\frac{1}{(z_{n,j}-z_{n-1,k})^2}}- \sum_{k=1,\; k\neq j}^{\frac{1}{2}n(n+1)}{\frac{1}{(z_{n,j}-z_{n,k})^2}}=-\frac{z_{n,j}}{6},\\
 \sum_{k=1}^{\frac{1}{2}n(n-1)}{\frac{1}{(z_{n,j}-z_{n-1,k})^3}}- \sum_{k=1,\; k\neq j}^{\frac{1}{2}n(n+1)}{\frac{1}{(z_{n,j}-z_{n,k})^3}}=-\frac{n-1}{4},\\
 \sum_{k=1}^{\frac{1}{2}n(n-1)}{\frac{1}{(z_{n,j}-z_{n-1,k})^5}}- \sum_{k=1,\; k\neq j}^{\frac{1}{2}n(n+1)}{\frac{1}{(z_{n,j}-z_{n,k})^5}}=z_{n,j}\left(\frac{n-1}{24}+\frac{1}{36}\right).
 \end{gather*}
\end{theorem}

\begin{proof}
Let $1\!\leq\! j\!\leq \!\frac{1}{2}n(n-1)$ and def\/ine $\omega:=z_{n-1,j}$ and $u:=w_n\!-\!\frac{1}{z-\omega}$. Since \mbox{$\gcd({Q_{n-1},Q_n})=1$}, see  Fukutani, Okamoto, and Umemura~\cite{fukutani}, equation~\eqref{eqyn} shows that $u$ is holomorphic in a~neighbourhood of $\omega$. Hence $u$ has a power series expansion, say
\[ \sum_{m=0}^{\infty}{a_m(z-\omega)^m},\]
which converges in an open disc centered at $\omega$.

Since $w_n$ is a solution of $P_{\rm II}(n)$, $u$ satisf\/ies
\begin{gather*}
(z-\omega)^2 u''=6u+6(z-\omega)u^2+2(z-\omega)^2u^3+(n+1)(z-\omega)^2+\omega(z-\omega)\\
\hphantom{(z-\omega)^2 u''=}{} +(z-\omega)^3u+\omega(z-\omega)^2u.
\end{gather*}
Hence we have the following identity in an open disc centered at $\omega$:
\begin{gather*}
 \sum_{m=2}^{\infty}{(m-1)ma_m(z-\omega)^m}=  6 \sum_{m=0}^{\infty}{a_m(z-\omega)^m}+6(z-\omega)\left( \sum_{m=0}^{\infty}{a_m(z-\omega)^m}\right)^2\\
\qquad{} +2(z-\omega)^2\left( \sum_{m=0}^{\infty}{a_m(z-\omega)^m}\right)^3+(n+1)(z-\omega)^2+\omega(z-\omega)\\
\qquad{} +(z-\omega)^3 \sum_{m=0}^{\infty}{a_m(z-\omega)^m}+\omega(z-\omega)^2 \sum_{m=0}^{\infty}{a_m(z-\omega)^m}.
\end{gather*}
By considering coef\/f\/icients of $(z-\omega)^n$, $n=0,1,2,4$, it is easy to deduce that $a_0=0$, $a_1=-\frac{\omega}{6}$, $a_2=-\frac{n+1}{4}$ and $a_4=\omega\left(\frac{n+1}{24}-\frac{1}{36}\right)$. Note that $a_3$ does not follow from considering coef\/f\/icients of $(z-\omega)^3$.

By Taylor's theorem and equation~\eqref{eqyn},
\begin{gather*}
a_m=\frac{u^{(m)}(z_{n-1,j})}{m!}=(-1)^m\left(  \sum_{k=1,\; k\neq j}^{\frac{1}{2}n(n-1)}\!\! {\frac{1}{(z_{n-1,j}-z_{n-1,k})^{m+1}}}- \! \sum_{k=1}^{\frac{1}{2}n(n+1)}\!\! {\frac{1}{(z_{n-1,j}-z_{n,k})^{m+1}}}\right)\!.
\end{gather*}
The f\/irst half of the theorem follows, the second half is proved analogously.
\end{proof}

Note that countably many nontrivial relations can be found between the $a_m$ in the above proof, by considering the coef\/f\/icient of $(z-\omega)^n$, for $n\in\mathbb{N}$.

In Kudryashov and Demina \cite{kudryashov} similar relations for the roots of $Q_n$ are obtained using the Korteweg--de Vries equation. In particular, the following results are presented in \cite{kudryashov} for $1\leq j \leq \frac{1}{2}n(n+1)$:
\begin{gather*}
\sum_{k=1,\; k\neq j}^{\frac{1}{2}n(n+1)}{\frac{1}{(z_{n,j}-z_{n,k})^2}}=-\frac{z_{n,j}}{12},\qquad
 \sum_{k=1,\; k\neq j}^{\frac{1}{2}n(n+1)}{\frac{1}{(z_{n,j}-z_{n,k})^3}}=0,\\
 \sum_{k=1,\; k\neq j}^{\frac{1}{2}n(n+1)}{\frac{1}{(z_{n,j}-z_{n,k})^5}}=-\frac{z_{n,j}}{144}.
 \end{gather*}
From these relations and Theorem \ref{TREL}, we obtain the following corollary:
\begin{corollary} For $1\leq j\leq \frac{1}{2}n(n-1)$:
\begin{gather*}
\sum_{k=1}^{\frac{1}{2}n(n+1)}{\frac{1}{(z_{n-1,j}-z_{n,k})^2}}=-\frac{z_{n-1,j}}{4},\qquad
 \sum_{k=1}^{\frac{1}{2}n(n+1)}{\frac{1}{(z_{n-1,j}-z_{n,k})^3}}=\frac{n+1}{4},\\
 \sum_{k=1}^{\frac{1}{2}n(n+1)}{\frac{1}{(z_{n-1,j}-z_{n,k})^5}}=-z_{n-1,j}\left(\frac{n+1}{24}-\frac{1}{48}\right).
 \end{gather*}
For $1\leq j\leq \frac{1}{2}n(n+1)$:
\begin{gather*}
\sum_{k=1}^{\frac{1}{2}n(n-1)}{\frac{1}{(z_{n,j}-z_{n-1,k})^2}}=-\frac{z_{n,j}}{4},\qquad
 \sum_{k=1}^{\frac{1}{2}n(n-1)}{\frac{1}{(z_{n,j}-z_{n-1,k})^3}}=-\frac{n-1}{4},\\
 \sum_{k=1}^{\frac{1}{2}n(n-1)}{\frac{1}{(z_{n,j}-z_{n-1,k})^5}}=z_{n,j}\left(\frac{n-1}{24}+\frac{1}{48}\right).
 \end{gather*}
\end{corollary}

In Theorem \ref{TREL}, we have obtained $4$ times $\frac{1}{2}n(n-1)$ plus $4$ times $\frac{1}{2}n(n+1)$ equations satisf\/ied by the $\frac{1}{2}n(n+1)$ roots of~$Q_n$, suggesting that these equations can be used to determine the roots of the polynomials~$Q_n$ recursively. If so, then these equations may be of use to derive properties of the roots of the Yablonskii--Vorob'ev polynomials. We shall not pursue this issue further here.

\section{Sums of negative powers of roots}\label{section4}

In Section~\ref{section2}, the rational solutions $w_n$ of $P_{\rm II}(n)$ were studied around roots of the Yablonskii--Vorob'ev polynomials. In this section, we consider $w_n$ at $0$.

Let $n\equiv 0 \pmod{3}$, then $0$ is not a root of $Q_{n-1}$ or $Q_n$. Therefore, by equation \eqref{eqyn}, $w_n$ is holomorphic in a neighbourhood of $0$. So $w_n$ has a power series expansion, say \[ \sum_{m=0}^{\infty}{a_mz^m},\]
which converges on an open disc centered at~$0$.

By Taylor's theorem and equation~\eqref{eqyn}, we have
\[a_m=-\left( \sum_{k=1}^{\frac{1}{2}n(n-1)}{\frac{1}{z_{n-1,k}^{m+1}}}- \sum_{k=1}^{\frac{1}{2}n(n+1)}{\frac{1}{z_{n,k}^{m+1}}}\right).\]
Let $\omega:=e^{\frac{2\pi i}{3}}$. Since $n\equiv 0 \pmod{3}$, $Q_n\in\mathbb{Z}[z^3]$. Therefore, the roots of $Q_n$ are invariant under multiplication by $\omega$. Hence
\[ \sum_{k=1}^{\frac{1}{2}n(n+1)}{\frac{1}{z_{n,k}^{m+1}}}=
\sum_{k=1}^{\frac{1}{2}n(n+1)}{\frac{1}{(\omega z_{n,k})^{m+1}}}=\frac{1}{\omega^{m+1}}
\sum_{k=1}^{\frac{1}{2}n(n+1)}{\frac{1}{z_{n,k}^{m+1}}},\]
therefore, if $m\not\equiv 2 \pmod{3}$,
\begin{gather}
\label{eqzero}
 \sum_{k=1}^{\frac{1}{2}n(n+1)}{\frac{1}{z_{n,k}^{m+1}}}=0.
\end{gather}
By the same reason, if $m\not\equiv 2 \pmod{3}$,
\[
\sum_{k=1}^{\frac{1}{2}n(n-1)}{\frac{1}{z_{n-1,k}^{m+1}}}=0.
\]
So $a_m=0$, if $m\not\equiv 2 \pmod{3}$, and in an open disc centered at $0$,
\[w_n(z)=  \sum_{m=0}^{\infty}{a_{3m+2}z^{3m+2}}.\]

Since $w_n$ is a solution of $P_{\rm II}(n)$, we have the following identity in an open disc centered at $0$:
\[
\sum_{m=0}^{\infty}{(3m+1)(3m+2)a_{3m+2}z^{3m}}=2\left(  \sum_{m=0}^{\infty}{a_{3m+2}z^{3m+2}}\right)^3+ \sum_{m=0}^{\infty}{a_{3m+2}z^{3m+3}}+n.
\]
Comparing coef\/f\/icients gives $a_2=\frac{1}{2}n$, $a_5=\frac{1}{40}n$ and $a_8=\frac{1}{2240}n+\frac{1}{224}n^3$.
We have obtained the following relations for $n\equiv 0 \pmod{3}$:
\begin{gather*}
\sum_{k=1}^{\frac{1}{2}n(n-1)}{\frac{1}{z_{n-1,k}^3}}- \sum_{k=1}^{\frac{1}{2}n(n+1)}{\frac{1}{z_{n,k}^3}}=-\frac{n}{2},\qquad
\sum_{k=1}^{\frac{1}{2}n(n-1)}{\frac{1}{z_{n-1,k}^6}}-\sum_{k=1}^{\frac{1}{2}n(n+1)}{\frac{1}{z_{n,k}^6}}=-\frac{n}{40},\\
 \sum_{k=1}^{\frac{1}{2}n(n-1)}{\frac{1}{z_{n-1,k}^9}}- \sum_{k=1}^{\frac{1}{2}n(n+1)}{\frac{1}{z_{n,k}^9}}=-\frac{1}{2240}n-\frac{1}{224}n^3.
\end{gather*}

If $n\equiv 1 \pmod{3}$, then $u:=w_n+\frac{1}{z}$ is holomorphic at $0$ and satisf\/ies
\[
z^2u''=6u-6zu^2+2z^2u^3+z^3u+(n-1)z^2.
\]
By considering the power series expansion of $u=w_n+\frac{1}{z}$ around $0$, the following relations are found:
\begin{gather*}
\sum_{k=1}^{\frac{1}{2}n(n-1)}{\frac{1}{z_{n-1,k}^3}}-
\sum_{k=1,\; z_{n,k}\neq 0 }^{\frac{1}{2}n(n+1)}{\frac{1}{z_{n,k}^3}}=\frac{1}{4}(n-1),\\
\sum_{k=1}^{\frac{1}{2}n(n-1)}{\frac{1}{z_{n-1,k}^6}}- \sum_{k=1,\; z_{n,k}\neq 0}^{\frac{1}{2}n(n+1)}{\frac{1}{z_{n,k}^6}}=\frac{1}{56}(n-1)+\frac{3}{112}(n-1)^2,\\
\sum_{k=1}^{\frac{1}{2}n(n-1)}{\frac{1}{z_{n-1,k}^9}}- \sum_{k=1,\; z_{n,k}\neq 0}^{\frac{1}{2}n(n+1)}{\frac{1}{z_{n,k}^9}}=\frac{1}{2800}(n-1)+\frac{9}{5600}(n-1)^2+\frac{1}{448}(n-1)^3.
\end{gather*}

If $n\equiv 2 \pmod{3}$, then $u:=w_n-\frac{1}{z}$ is holomorphic at $0$ and satisf\/ies
\[
z^2u''=6u-6zu^2+2z^2u^3+z^3u+(n+1)z^2.
\]
By considering the power series expansion of $u=w_n-\frac{1}{z}$ around $0$, the following relations are found:
\begin{gather*}
\sum_{k=1,\; z_{n-1,k}\neq 0}^{\frac{1}{2}n(n-1)}{\frac{1}{z_{n-1,k}^3}}- \sum_{k=1}^{\frac{1}{2}n(n+1)}{\frac{1}{z_{n,k}^3}}=\frac{1}{4}(n+1),\\
 \sum_{k=1,\; z_{n-1,k}\neq 0}^{\frac{1}{2}n(n-1)}{\frac{1}{z_{n-1,k}^6}}- \sum_{k=1}^{\frac{1}{2}n(n+1)}{\frac{1}{z_{n,k}^6}}=\frac{1}{56}(n+1)-\frac{3}{112}(n+1)^2,\\
 \sum_{k=1,\; z_{n-1,k}\neq 0}^{\frac{1}{2}n(n-1)}{\frac{1}{z_{n-1,k}^9}}- \sum_{k=1}^{\frac{1}{2}n(n+1)}{\frac{1}{z_{n,k}^9}}=\frac{1}{2800}(n+1)-\frac{9}{5600}(n+1)^2+\frac{1}{448}(n+1)^3.
\end{gather*}

\begin{remark}
\label{remark}
Considering higher order coef\/f\/icients, we see that for every threefold $m\geq 3$, polynomial expressions in $n$, with rational coef\/f\/icients, depending on $n \pmod{3}$, exist for
\[ \sum_{k=1, \; z_{n-1,k}\neq 0}^{\frac{1}{2}n(n-1)}{\frac{1}{z_{n-1,k}^m}}- \sum_{k=1, \; z_{n,k}\neq 0}^{\frac{1}{2}n(n+1)}{\frac{1}{z_{n,k}^m}}.\]
\end{remark}

As a corollary of these relations, by induction, we obtain:
\begin{gather*}
\sum_{k=1,\; z_{n,k}\neq 0 }^{\frac{1}{2}n(n+1)}{\frac{1}{z_{n,k}^3}}=\begin{cases} \dfrac{n}{4}& \text{if $n\equiv 0 \pmod{3}$,}\vspace{2mm}\\
0& \text{if $n\equiv 1 \pmod{3}$,}\vspace{2mm}\\
-\dfrac{n+1}{4}& \text{if $n\equiv 2 \pmod{3},$}
\end{cases}
\\
\sum_{k=1,\;z_{n,k}\neq 0 }^{\frac{1}{2}n(n+1)}{\frac{1}{z_{n,k}^6}}=\begin{cases} \dfrac{1}{40}n^2+\dfrac{1}{80}n& \text{if $n\equiv 0 \pmod{3}$,}\vspace{2mm}\\
-\dfrac{1}{560}n^2-\dfrac{1}{560}n+\dfrac{1}{280}& \text{if $n\equiv 1 \pmod{3}$,}\vspace{2mm}\\
\dfrac{1}{40}n^2+\dfrac{3}{80}n+\dfrac{1}{80}& \text{if $n\equiv 2 \pmod{3}$,}
\end{cases}
\\
\sum_{k=1,\; z_{n,k}\neq 0 }^{\frac{1}{2}n(n+1)}{\frac{1}{z_{n,k}^9}}=\begin{cases} \dfrac{n+7n^2+10n^3}{4480}& \text{if $n\equiv 0 \pmod{3}$,}\vspace{2mm}\\
\dfrac{2-n-n^2}{22400}& \text{if $n\equiv 1 \pmod{3}$,}\vspace{2mm}\\
\dfrac{-20-85n-115n^2-50n^3}{22400}& \text{if $n\equiv 2 \pmod{3}$.}
\end{cases}
\end{gather*}
By Remark \ref{remark}, for every threefold $m\geq 3$, polynomial expressions in $n$, with rational coef\/f\/icients, depending on $n \pmod{3}$, exist for
\[
\sum_{k=1,\; z_{n,k}\neq 0 }^{\frac{1}{2}n(n+1)}{\frac{1}{z_{n,k}^m}}.
\]
If $m\not\equiv 0 \pmod{3}$, see equation~\eqref{eqzero}, then
\[
\sum_{k=1,\;z_{n,k}\neq 0 }^{\frac{1}{2}n(n+1)}{\frac{1}{z_{n,k}^m}}=0.
\]
So, for all $n,m\in\mathbb{N}$, \[
\sum_{k=1,\; z_{n,k}\neq 0 }^{\frac{1}{2}n(n+1)}{\frac{1}{z_{n,k}^m}}\in\mathbb{Q},
\]
even though the nonzero roots of the Yablonskii--Vorob'ev polynomials are irrational.

\subsection*{Acknowledgements}
I wish to thank Erik Koelink for his enlightening discussions and introducing me to the world of the Painlev\'e equations. I am also grateful to Peter Clarkson for his interest and useful links to the literature.

\pdfbookmark[1]{References}{ref}
\LastPageEnding

\end{document}